\documentstyle[txmac,a4,
amssymb,%
case,%
twoside,%
nocaphead,%
epsf,%
varthm,%
myrot,%
mypic,times,mathptm]{article}

\advance\oddsidemargin by -1.9cm
\advance\evensidemargin by -1.9cm
\advance\textwidth by 3.8cm

\def\mynewtheo#1#2{%
\newtheorem{@#1}{#2}[section]%
\newenvironment{#1}{\begin{@#1}\rm}{\end{@#1}}}

\mynewtheo{lemma}{Lemma}
\mynewtheo{exer}{Exercise}
\mynewtheo{theo}{Theorem}
\mynewtheo{remark}{Remark}
\mynewtheo{defi}{Definition}
\mynewtheo{conj}{Conjecture}
\mynewtheo{corollary}{Corollary}
\mynewtheo{prop}{Proposition}
\mynewtheo{question}{Question}
\mynewtheo{exam}{Example}

\newenvironment{rem}{\begin{remark}}{\end{remark}}
\newenvironment{theorem}{\begin{theo}}{\end{theo}}
\newenvironment{conjecture}{\begin{conj}}{\end{conj}}

\begin{document}

\makeatletter

\newenvironment{eqn}{\begin{equation}}{\end{equation}\@ignoretrue}

\newenvironment{myeqn*}[1]{\begingroup\def\@eqnnum{\reset@font\rm#1}%
\xdef\@tempk{\arabic{equation}}\begin{equation}\edef\@currentlabel{#1}}
{\end{equation}\endgroup\setcounter{equation}{\@tempk}\ignorespaces}

\newenvironment{myeqn}[1]{\begingroup\let\eq@num\@eqnnum
\def\@eqnnum{\bgroup\let\r@fn\normalcolor 
\def\normalcolor####1(####2){\r@fn####1#1}%
\eq@num\egroup}%
\xdef\@tempk{\arabic{equation}}\begin{equation}\edef\@currentlabel{#1}}
{\end{equation}\endgroup\setcounter{equation}{\@tempk}\ignorespaces}

\newenvironment{myeqn**}{\begin{myeqn}{(\arabic{equation})\es\es\mbox{\qed}}\edef\@currentlabel{\arabic{equation}}}
{\end{myeqn}\stepcounter{equation}}

\def\myfrac#1#2{\raisebox{0.2em}{\small$#1$}\!/\!\raisebox{-0.2em}{\small$#2$}}

\def\myeqnlabel{\bgroup\@ifnextchar[{\@maketheeq}{\immediate
\stepcounter{equation}\@myeqnlabel}}

\def\@maketheeq[#1]{\def\theequation{#1}\@myeqnlabel}

\def\@myeqnlabel#1{%
{\edef\@currentlabel{\theequation}
\label{#1}\enspace\eqref{#1}}\egroup}

\def\epsfs#1#2{{\ifautoepsf\unitxsize#1\relax\else
\epsfxsize#1\relax\fi\epsffile{#2.eps}}}

\def\@test#1#2#3#4{%
  \let\@tempa\go@
  \@tempdima#1\relax\@tempdimb#3\@tempdima\relax\@tempdima#4\unitxsize\relax
  \ifdim \@tempdimb>\z@\relax
    \ifdim \@tempdimb<#2%
      \def\@tempa{\@test{#1}{#2}}%
    \fi
  \fi
  \@tempa
}

\def\bysame{\same[\kern2cm]\,}

\def\go@#1\@end{}
\newdimen\unitxsize
\newif\ifautoepsf\autoepsftrue

\unitxsize4cm\relax
\def\epsfsize#1#2{\epsfxsize\relax\ifautoepsf
  {\@test{#1}{#2}{0.1 }{4   }
		{0.2 }{3   }
		{0.3 }{2   }
		{0.4 }{1.7 }
		{0.5 }{1.5 }
		{0.6 }{1.4 }
		{0.7 }{1.3 }
		{0.8 }{1.2 }
		{0.9 }{1.1 }
		{1.1 }{1.  }
		{1.2 }{0.9 }
		{1.4 }{0.8 }
		{1.6 }{0.75}
		{2.  }{0.7 }
		{2.25}{0.6 }
		{3   }{0.55}
		{5   }{0.5 }
		{10  }{0.33}
		{-1  }{0.25}\@end
		\ea}\ea\epsfxsize\the\@tempdima\relax
		\fi
		}


\let\old@tl\~
\def\~{\raisebox{-0.8ex}{\tt\old@tl{}}}

\author{A. Stoimenow\footnotemark[1]\\[2mm]
\small Department of Mathematics, \\
\small University of Toronto,\\
\small Canada M5S 3G3\\
\small e-mail: {\tt stoimeno@math.toronto.edu}\\
\small WWW: {\hbox{\tt http://www.math.toronto.edu/stoimeno/}}
}

{\def\thefootnote{\fnsymbol{footnote}}
\footnotetext[1]{Supported by a DFG postdoc grant.}
}

\title{\large\bf \uppercase{Polynomial values, the linking form and
unknotting numbers}\\[4mm]
{\small\it This is a preprint. I would be grateful for any comments
and corrections!}
}

\date{\phantom{\large Current version: \curv\ \ \ First version:
\makedate{1}{2}{2000}}}

\maketitle

\makeatletter

\let\point\pt
\let\ay\asymp
\let\pa\partial
\let\al\alpha
\let\be\beta
\let\Gm\Gamma
\let\gm\gamma
\let\de\delta
\let\Dl\Delta
\let\eps\epsilon
\let\lm\lambda
\let\Lm\Lambda
\let\sg\sigma
\let\vp\varphi
\let\zt\zeta
\let\om\omega

\let\reference\ref
\let\lm\lambda
\let\sg\sigma
\def\bQ{{\Bbb Q}}
\def\bN{{\Bbb N}}
\def\bZ{{\Bbb Z}}

\let\dt\det
\let\sm\setminus
\let\tl\tilde

\def\ncap{\not\mathrel{\cap}}
\def\cf{\text{\rm cf}\,}
\def\lra{\longrightarrow}
\def\so{\Rightarrow}
\def\So{\Longrightarrow}
\def\llra{\longleftrightarrow}
\let\ds\displaystyle
\def\bt{\bar t_2}
\def\TM{$^\text{\raisebox{-0.2em}{${}^\text{TM}$}}$}
\def\ppmod{\kern-0.85em\pmod}

\let\reference\ref

\long\def\@makecaption#1#2{%
   \vskip 10pt
   {\let\label\@gobble
   \let\ignorespaces\@empty
   \xdef\@tempt{#2}%
   }%
   \ea\@ifempty\ea{\@tempt}{%
   \setbox\@tempboxa\hbox{%
      \fignr#1#2}%
      }{%
   \setbox\@tempboxa\hbox{%
      {\fignr#1:}\capt\ #2}%
      }%
   \ifdim \wd\@tempboxa >\captionwidth {%
      \rightskip=\@captionmargin\leftskip=\@captionmargin
      \unhbox\@tempboxa\par}%
   \else
      \hbox to\captionwidth{\hfil\box\@tempboxa\hfil}%
   \fi}%
\def\fignr{\small\sffamily\bfseries}%
\def\capt{\small\sffamily}%

\newdimen\@captionmargin\@captionmargin2cm\relax
\newdimen\captionwidth\captionwidth\hsize\relax

\def\eqref#1{(\protect\ref{#1})}

\def\proof{\@ifnextchar[{\@proof}{\@proof[\unskip]}}
\def\@proof[#1]{\noindent{\bf Proof #1.}\enspace}

\def\hint{\noindent Hint: }
\def\problem{\noindent{\bf Problem.} }

\def\@mt#1{\ifmmode#1\else$#1$\fi}
\def\qed{\hfill\@mt{\Box}}
\def\qqed{\hfill\@mt{\Box\enspace\Box}}

\def\cU{{\cal U}}
\def\cC{{\cal C}}
\def\cP{{\cal P}}
\def\tP{{\tilde P}}
\def\tZ{{\tilde Z}}
\def\fg{{\frak g}}
\def\sgn{\text{\operator@font sgn}}
\def\cZ{{\cal Z}}
\def\cD{{\cal D}}
\def\bR{{\Bbb R}}
\def\cE{{\cal E}}
\def\bQ{{\Bbb Q}}
\def\bZ{{\Bbb Z}}
\def\bN{{\Bbb N}}

\def\abstractname{}

\@addtoreset {footnote}{page}

\renewcommand{\section}{%
   \@startsection
         {section}{1}{\z@}{-1.5ex \@plus -1ex \@minus -.2ex}%
               {1ex \@plus.2ex}{\large\bf}%
}
\renewcommand{\@seccntformat}[1]{\csname the#1\endcsname .
\quad}

\def\bC{{\Bbb C}}
\def\bP{{\Bbb P}}

{\let\@noitemerr\relax
\vskip-2.7em\kern0pt\begin{abstract}
\noindent{\bf Abstract.}\enspace
We show how the signed evaluations of link polynomials can be used to
calculate unknotting numbers. We use the 
Jones-Rong value of the Brandt-Lickorish-Millett-Ho polynomial $Q$ to
calculate the unknotting numbers of $8_{16}$, $9_{49}$ and
6 further new entries in Kawauchi's tables. Another method is
developed by applying and extending the linking form criterion of
Lickorish. This leads to several conjectured
relations between the Jones value of $Q$ and the linking form.\\[1mm]

(secondary). \\[1mm]
\end{abstract}
}

\def\epsfsv#1#2{{\vcbox{\epsfs{#1}{#2}}}}
\def\vcbox#1{\setbox\@tempboxa=\hbox{#1}\parbox{\wd\@tempboxa}{\box
  \@tempboxa}}

\def\fr#1{\left\lfloor#1\right\rfloor}
\def\br#1{\left\langle#1\right\rangle}
\def\BR#1{\left\lceil#1\right\rceil}

\section{Introduction}

The unknotting number $u(K)$ of a knot $K$ is defined as the minimal
number of crossing changes in any diagram of $K$ needed to turn
$K$ into the unknot (see e.g. \cite{ACampo,Kawamura,Kobayashi,%
Nakanishi,Scharl,Tanaka,Wendt,Zhang}). Despite its simple definition,
the unknotting number has proved hard to calculate. The tables in
\cite{Kawauchi} show that, after several decades of work, the
unknotting number of each fifth prime knot with 10 or fewer crossings
remains unknown.

In this paper, we develop two at first glance different approaches
to the calculation of unknotting numbers~-- using the
evaluations of the link polynomials \cite{Jones,LickMil,Lipson,Rong}
and the linking form $\lm$ on $H_1(D_K)$ \cite{Lickorish}.

As outcome, we settle the problem of determining the unknotting
number for 9 (that is, about one sixth) of the open nine and ten
crossing knots in Kawauchi's tables \cite{Kawauchi}.
The table below summarizes these examples:
\begin{eqn}\label{teq}
\begin{mytab}{c||l|l|l}{ & \multicolumn{3}{|c}{ } }
method &
\multicolumn{1}{|c|}{$u=1$} &
\multicolumn{1}{|c}{$u=2$} &
\multicolumn{1}{|c}{$u=3$}
\\[2mm]%
\hline
\hline[2mm]%
observation & $10_{131}$ & & \\
$Q$ &  & $8_{16}$, $10_{86}$, \phantom{$10_{105}$, }$10_{106}$, $10_{109}$, $10_{116}$,
$10_{121}$ & $9_{49}$, $10_{103}$ \\[2mm]%
$\lm$ & & $8_{16}$, $10_{86}$, $10_{105}$, $10_{106}$, $10_{109}$,
$10_{116}$, $10_{121}$ & $9_{49}$ \\[2mm]%
\end{mytab}
\end{eqn}

Our criteria can also be applied to some simple composite knots
(we give a table of such unknotting numbers in an appendix) and to
knot distance \cite{Murakami}.

It is striking that the outcome of both methods~-- the $Q$
polynomial and the linking form~-- give surprisingly similar, although
not identical, results. This is clearly hardly a matter of
accidental coincidence, and thus we are led to several conjectures
on relations between both. We will mention the evidence for most of
these conjectures while discussing the various examples, and explicitly
compile the conjectures for the benefit of the reader at the end of the
paper, hoping to motivate further investigations on this subject.

\section{Preliminaries and notation}

For each knot $K$ we have a sequence of knots $K_i$
\begin{eqn}\label{3K}
K_0\to K_1\to K_2\to\dots\to K_n
\end{eqn}
such that $K=K_0$, $K_n$ is the unknot, and $K_i$ differs from
$K_{i-1}$ only by a crossing change.
We call a sequence \eqref{3K} an \em{unknotting sequence} for $K$.
The minimal length $n$ of an unknotting sequence \eqref{3K} for $K$
is the \em{unknotting number} $u(K)$ of $K$.

Henceforth, $D_K$ denotes the \em{double branched cover} of
$S^3$ over a knot $K$. By $H_1=H_1(D_K)=H_1(D_K,\bZ)$
we denote its homology group over $\bZ$. (The various abbreviated
versions will be used at places where no confusion arises; $H_1$
will be used throughout the paper only in this context, so that,
for example, when we talk of $H_1$ of a knot, always $H_1$ of its
double cover will be meant.) $H_1$ is a finite commutative group of
odd order. This order is called the \em{determinant} of a knot
$K$, and it will be denoted as $\dt(K)$. (This generalizes to
links $L$, by putting $\dt(L)=0$ to stand for infinite $H_1(D_L)$.)
By the classification of finite commutative groups, $H_1$ decomposes
into a direct sum of finite (odd order) cyclic groups $\bZ_p$; their
orders are called \em{torsion numbers}. Wendt \cite{Wendt} proved
that the number of torsion numbers of $H_1(D_K)$ is not smaller 
than $u(K)$. $H_1$ is also equipped with a bilinear form
$\lm\,:\,H_1\times H_1 \to\bQ/\bZ$, called the \em{linking form}
(see \cite{Lickorish,MurYas} for example).

In the following knots and links will be assumed oriented, but
sometimes orientation will be irrelevant.

The \em{Jones polynomial} $V$ (introduced in \cite{Jones}, but now
commonly used with the convention of \cite{Jones3})
is a Laurent polynomial in one
variable $t$ of oriented knots and links, and can be defined
by being $1$ on the unknot and the \em{(skein) relation}
\begin{eqn}\label{1} 
t^{-1}\,V\left(
L_+
\right)\,-\,
t \,V\left(
L_-
\right)\,=\,
-(t^{-1/2}-t^{1/2})\,V\left(
L_0
\right)\,.
\end{eqn}
Herein $L_{\pm,0}$ are three links with diagrams differing only near
a crossing. 
\begin{eqn}\label{Lpm0}
\begin{array}{*2{c@{\qquad}}c}
\diag{9mm}{1}{1}{
\picmultivecline{-8 1 -1.0 0}{1 0}{0 1}
\picmultivecline{-8 1 -1.0 0}{0 0}{1 1}
}
 &
\diag{9mm}{1}{1}{
\picmultivecline{-8 1 -1 0}{0 0}{1 1}
\picmultivecline{-8 1 -1 0}{1 0}{0 1}
}
&
\diag{9mm}{1}{1}{
\piccirclevecarc{1.35 0.5}{0.7}{-230 -130}
\piccirclevecarc{-0.35 0.5}{0.7}{310 50}
}
\\[2mm]
L_+ & L_- & L_0
\end{array}
\end{eqn}
We call the crossings in the first two fragments resp.\ \em{positive}
and \em{negative}, and a crossing replaced by the third fragment
\em{smoothed out}. A triple of links that can be represented as
$L_{\pm,0}$ in \eqref{Lpm0} is called a \em{skein triple}. The sum
of the signs ($\pm 1$) of the crossings of a diagram $D$
is called \em{writhe} of $D$ and written $w(D)$.

\begin{figure}[htb]
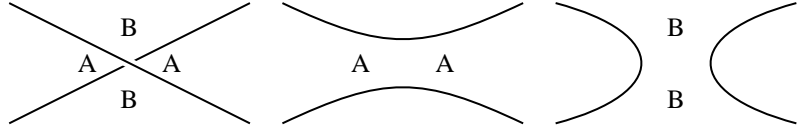

\[
\diag{8mm}{4}{2}{
   \picline{0 0}{4 2}
   \picmultiline{-5.0 1 -1.0 0}{0 2}{4 0}
   \picputtext{2.7 1}{A}
   \picputtext{1.3 1}{A}
   \picputtext{2 1.6}{B}
   \picputtext{2 0.4}{B}
} \quad
\diag{8mm}{4}{2}{
   \pictranslate{2 1}{
       \picmultigraphics[S]{2}{1 -1}{
           \piccurve{-2 1}{-0.3 0.2}{0.3 0.2}{2 1}
       }
   }
   \picputtext{2.7 1}{A}
   \picputtext{1.3 1}{A}
} \quad
\diag{8mm}{4}{2}{
   \pictranslate{2 1}{
       \picmultigraphics[S]{2}{-1 1}{
           \piccurve{2 -1}{0.1 -0.5}{0.1 0.5}{2 1}
       }
   }
   \picputtext{2 1.6}{B}
   \picputtext{2 0.4}{B}
}
\]
\caption{\label{figsplit}The A- and B-corners of a
crossing, and its both splittings. The corner A (resp. B)
is the one passed by the overcrossing strand when rotated 
counterclockwise (resp. clockwise) towards the undercrossing 
strand. A type A (resp.\ B) splitting is obtained by connecting 
the A (resp.\ B) corners of the crossing.}
\end{figure}

A different interpretation of the Jones polynomial
than via skein rules has been developed by Kauffman
\cite{Kauffman} (see also \cite[\S 6.2]{Adams}). The Kauffman
state model is sometimes more useful than the skein
approach, and we shall also consider it below. Recall, that
the \em{Kauffman bracket} $\br{D}$ of a(n unoriented) link
diagram $D$ is a Laurent polynomial in a variable
$A$, obtained by summing over all states the terms
\begin{eqn}\label{eq_12}
A^{\#A-\#B}\,\left(-A^2-A^{-2}\right)^{|S|-1}\,.
\end{eqn}
Herein a \em{state} is a choice of \em{splittings} of type $A$ or 
$B$ for any single crossing (see figure \ref{figsplit}), 
$\#A$ and $\#B$ denote the number of
type A (resp. type B) splittings and $|S|$ the
number of (disjoint) circles obtained after all
splittings in a state.

The Jones polynomial of a link $L$ can be calculated from the
Kauffman bracket, by evaluating it on the (unoriented version of)
a diagram $D$ of $L$, and then multiplying by a power of $t$
coming from the (orientation dependent) writhe of $D$:
\begin{eqn}\label{conv}
V_L(t)\,=\,\left(-t^{-3/4}\right)^{-w(D)}\,\br{D}
\raisebox{-0.6em}{$\Big |_{A=t^{-1/4}}$}\,.
\end{eqn}

The \em{signature} $\sg$ is a $\bZ$-valued invariant of knots and links.
Originally it was defined terms of Seifert matrices \cite{Rolfsen}.
We have that $\sg(L)$ has the opposite parity to the number of
components of a link $L$, whenever the determinant of $L$ is non-zero
(i.e. $H_1(D_L)$ is finite). This in particular always happens
for $L$ being a knot, so that $\sg$ takes only even values on knots.

The most of the early work on the signature was done by
Murasugi \cite{Murasugi}, who showed several properties
of this invariant. In particular the following property
is known: if $L_{\pm,0}$ form a skein triple, then
\begin{eqnarray}
\sg(L_+)-\sg(L_-) & \in & \{0,1,2\}\,, \label{2a} \\[2mm]
\sg(L_\pm)-\sg(L_0) & \in & \{-1,0,1\}\,. \label{2b}
\end{eqnarray}
(Note: In \eqref{2a} one can also have $\{0,-1,-2\}$ instead
of $\{0,1,2\}$, since other authors, like Murasugi, take $\sg$ to
be with opposite sign. Thus \eqref{2a} not only defines a property,
but also specifies our sign convention for $\sg$.)
We remark that for knots in \eqref{2a} only $0$ and $2$ can occur
on the right. A consequence of this relation is the inequality
$u(K)\ge|\sg(K)/2|$. After Wendt's inequality, this was one of the
first important results on the unknotting number.

%
%

Now we introduce the \em{Brandt-Lickorish-Millett-Ho polynomial}
$Q_K=Q_K(z)$ of $K$ \cite{BLM,Ho}.
Recall, that the $Q$ polynomial is a Laurent polynomial in one
variable $z$ for links without orientation, defined by being 1 on the
unknot and the relation
\begin{eqn}\label{Qrel}
z\,\Bigl[\ 
Q\Bigl(\diag{3mm}{2}{2}{
    \piccirclearc{-1 1}{1.41}{-45 45}
    \piccirclearc{3 1}{1.41}{135 -135}
  }
\Bigr)+
Q\Bigl(\diag{3mm}{2}{2}{
  \piccirclearc{1 -1}{1.41}{45 135}
    \piccirclearc{1 3}{1.41}{-135 -45}
}\Bigr)\ \Bigr]\quad=\quad
Q\Bigl(\diag{3mm}{2}{2}{
	\picmultiline{-10.5 1 -1 0}{0 0}{2 2}
	\picmultiline{-10.5 1 -1 0}{2 0}{0 2}
    }\Bigr)+
Q\Bigl(\diag{3mm}{2}{2}{
  \picmultiline{-10.5 1 -1 0}{2 0}{0 2}
  \picmultiline{-10.5 1 -1 0}{0 0}{2 2}
}\Bigr)\,,
\end{eqn}
where again the fragments denote link
diagrams equal except in the specified fragment.

If we modify the skein relation for $V$ by omitting the coefficients
$t^{\mp 1}$ of $L_{\pm}$ on the left of \eqref{1}, we obtain the
skein relation for another (and more classical) polynomial invariant,
the \em{Alexander polynomial} $\Dl(t)$ (see \cite{Rolfsen}).

All three polynomials allow to express the determinant of $K$, as
\[
\dt(K)=\big|\,\Dl_K(-1)\,\big|\,=\,
\big|\,V_K(-1)\,\big|\,=\,\sqrt{Q_K(2)}\,.
\]
There are further special values of the Jones and $Q$ polynomial, which
will be discussed in the following (see \cite[\S 12]{Jones3},
\cite{LickMil} and \cite{BLM} for more details).

By $k\% 2\in\{0,1\}$ we denote the parity of $k$ and by
$\fr{k}$ the greatest integer not exceeding $k$. `W.l.o.g.' will
abbreviate `without loss of generality'.
The notation for knots we use is this of Rolfsen \cite{Rolfsen}.

\section{The Jones polynomial\label{JP}}

We first start by a property of the Jones polynomial $V$ of an
unknotting number one knot, which slightly generalizes Traczyk's
criterion for the unknotting number one case. It is related to the
\em{signed} unknotting number (see \cite{Traczyk,CL}).

\begin{prop}\label{thV}
Let $K$ be an unknotting number one knot which can be unknotted by
switching a positive crossing to the negative, and set
\begin{eqn}\label{tlV}
\tl V\,=\,1-\frac{V_K-1}{t-1}\,.
\end{eqn}
Then there exists a knot $K'$ with $V_{K'}=t^{-\tl V'(1)}\tl V$.
\end{prop}

\proof Let $D$ be the diagram of $K$ unknotting by a change of a
positive crossing. (We call this crossing unknotting crossing.)
W.l.o.g. we can assume $D$ to have zero writhe (add kinks). Thus,
considering the Kauffman bracket version of $V$, we have
\begin{eqnarray*}
\br{\diag{4mm}{2}{2}{
  \picmultiline{-8.5 1 -1 0}{2 0}{0 2}
  \picmultiline{-8.5 1 -1 0}{0 0}{2 2}
}} & = & V_K(A^{-4}) \quad\mbox{and}\\
\br{\diag{4mm}{2}{2}{
  \picmultiline{-8.5 1 -1 0}{0 0}{2 2}
  \picmultiline{-8.5 1 -1 0}{2 0}{0 2}
}} & = & A^{-6}\,.
\end{eqnarray*}
Resolving the crossing according to the Kauffman bracket relation,
we get
\[
\left(\begin{array}{cc} A & A^{-1} \\ A^{-1} & A\end{array}\right)
\left(\begin{array}{c} 
\br{\diag{2mm}{2}{2}{
  \piccirclearc{-1 1}{1.41}{-45 45}
  \piccirclearc{3 1}{1.41}{135 -135}
}}\\[2mm]
\br{\diag{2mm}{2}{2}{
  \piccirclearc{1 -1}{1.41}{45 135}
  \piccirclearc{1 3}{1.41}{-135 -45}
}}\end{array}\right)\,=\,
\left(\begin{array}{c} 
V_K(A^{-4})\\
A^{-6}\end{array}\right)\,.
\]
The sign of the crossing switched shows that from both splicings,
$\diag{4mm}{2}{2}{
  \piccirclearc{1 -1}{1.41}{45 135}
  \piccirclearc{1 3}{1.41}{-135 -45}
}$ is the one that again corresponds to a knot. 
We find
\begin{eqnarray*}
\br{\diag{4mm}{2}{2}{
  \piccirclearc{1 -1}{1.41}{45 135}
  \piccirclearc{1 3}{1.41}{-135 -45}
}} & = & -\frac{V_K(A^{-4})A^{-1}-A\cdot A^{-6}}{A^2-A^{-2}} \\
& = & A^{-3}\left[\,\frac{-A^{-4}+V_K(A^{-4})}{A^{-4}-1}\,\right]\,,
\end{eqnarray*}
and setting $t=A^{-4}$ we get
\[
\br{\diag{4mm}{2}{2}{
  \piccirclearc{1 -1}{1.41}{45 135}
  \piccirclearc{1 3}{1.41}{-135 -45}
}} = -A^{-3}\left[1-\frac{V_K-1}{t-1}\,\right]\,.
\]
Thus the Jones polynomial of $K'=\diag{3mm}{2}{2}{
  \piccirclearc{1 -1}{1.41}{45 135}
  \piccirclearc{1 3}{1.41}{-135 -45}
}$ differs from $\tl V$ just by a unit in $Z[t,t^{-1}]$, and the
properties $V_{\tl K}(1)=1$ and $V'_{\tl K}(1)=0$ for $\tl K=K,K'$
show that this unit is precisely $t^{-\tl V'(1)}$. \qed

It has been known for a while, by \cite{LickMil}, that $V\left(
e^{\pi i/3}\right)=\pm (i\sqrt{3})^d$ for some $d\in\bN$,
which was observed by Traczyk in \cite{Traczyk} to show $u\ge d$
(a conclusion that alternatively follows when combining \cite{LickMil}
and \cite{Wendt}). He used a more intricate
combination of the sign of the evaluation and the signature to show
that sometimes $u>d$, proving $u(10_{67})=2$. Traczyk's observation
for $u=1$ is also contained in ours, but we can now apply more on what
we know on the values of $V$. We have the following properties,
basically due to Jones \cite[\S 12]{Jones3} and Lickorish-Millett
\cite{LickMil} (see also \cite{granny}).

\begin{prop}\label{th2}(Jones \cite[\S 12]{Jones3})
Let $V\in\bZ[t,t^{-1}]$ be that Jones polynomial of a knot. Then
$V\equiv 1\bmod(t-1)(t^3-1)$ and $V\equiv \pm 3^k(2t-1)^d \bmod t^2-t+1$
for some $k\in\bN$ and $d\in\{0,1\}$, such that if $3\nmid V(-1)$,
$k=d=0$ and if $3\mid V(-1)$, $k+d>0$ and $3^{2k+d}\mid V(-1)$.
Furthermore, we have the Arf invariant identity $V(i)=(-1)^{V''(1)/6}$.
\end{prop}

(We reformulated these conditions using the fact,
known from standard Galois theory, that
the evaluation of a polynomial $P$ in some algebraic number $v$ is
equivalent to the residue of $P$ modulo the minimal polynomial of $v$.)

%

A more general and simple special case is the following.

\begin{corollary}\label{cMi}
Let $K$ be an achiral knot (or weaker a knot with $V_K(t)=V_K(1/t)$),
and $3\mid V_K(-1)$. Then even $9\mid V_K(-1)$, and $u(K)>1$.
\end{corollary}

\proof If $K$ is achiral, then for $|t|=1$ we have $\bar t=-1/t$, and
thus $V_K(t)$ is real. If $u(K)=1$, then $\dim_{\bZ_3}H_1(D_K,\bZ_3)=1$,
but in this case $V_K\bigl(e^{\pi i/3}\bigr)=\pm i\sqrt{3}$, which is
not real. Therefore, $u(K)\ge \dim_{\bZ_3}H_1(D_K,\bZ_3)\ge 2$, so
in particular $9\mid V_K(-1)$. \qed

\begin{remark}
It was observed in \cite[\S 5]{granny}, that the property 
$\left(V_K\left(e^{\pi i/3}\right)\right)^2\mid V_K(-1)$ for the
Jones polynomial $V_K$ of a knot $K$ does not (in general) follow from
the properties of $V$ listed in \cite[\S 12]{Jones3}.
\end{remark}

A similar statement, which is unrelated to unknotting numbers, but
worth mentioning, is the following corollary. Here `$3$-equivalent'
means transformable by a sequence of Nakanishi's $3$-moves
(see \cite{granny}):
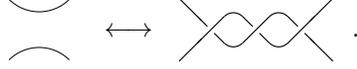
\[
\diag{4mm}{2}{2}{
  \piccirclearc{1 -1}{1.41}{45 135}
  \piccirclearc{1 3}{1.41}{-135 -45}
} 
\quad\llra\quad
\diag{6mm}{3.5}{1}{
  \piclinewidth{40}
  \picmultigraphics{3}{1 0}{
    \picline{0.2 0.8}{0.8 0.2}
    \picmultiline{-8.5 1 -1.0 0}{0.2 0.2}{0.8 0.8}
  }
  \picmultigraphics{2}{1 0}{
    \piccirclearc{1 0.6}{0.28}{45 135}
    \piccirclearc{1 0.4}{0.28}{-135 -45}
  }
  \picline{-0.2 -0.2}{0.3 0.3}
  \picline{-0.2 1.2}{0.3 0.7}
  \picline{3.2 -0.2}{2.7 0.3}
  \picline{3.2 1.2}{2.7 0.7}
}\,.
\]

\begin{corollary}
If a knot $K$ is $3$-equivalent to an unlink of an even number of
components, then $K$ is chiral. 
\end{corollary}

\proof If $K$ is $3$-equivalent to an unlink of an even number of
components, then $V_K\left(e^{\pi i/3}\right)$ is, up to a sign,
an odd power of $i\sqrt{3}$, which is not real. \qed



\section{The $Q$ polynomial\label{SQ}}

Now we go a step further from the proof of Proposition \reference{thV}
and consider the $Q$ polynomial \cite{BLM,Ho}.

\begin{prop}\label{th3}
Let $K$ be an unknotting number one knot of determinant
$\dt(K)=2n+1$, $n>0$. Set for $k\in\bN$
\[
S_5(k)\,:=\,\left\{\begin{array}{ll}\{0\} & 5\nmid k \\
\{l>0\,:\,5^l|k\} & 5\mid k \end{array}\right.
\]
Then there exist numbers $k\in S_5(n)$ and $l\in S_5(n+1)$ such that
\[
1+Q_K(z)\,\equiv\,z\left(\,\pm_15^{\fr{k/2}}(2z+1)^{k\% 2}
\pm_25^{\fr{l/2}}(2z+1)^{l\% 2}\right)\bmod z^2+z-1\,,
\]
for some (independent) sign choices $\pm_{1,2}\in\{+,-\}$.
\end{prop}

\proof We return to the arguments in the proof of Proposition
\reference{thV}. A similar calculation of $\br{\diag{3mm}{2}{2}{
  \piccirclearc{-1 1}{1.41}{-45 45}
  \piccirclearc{3 1}{1.41}{135 -135}
}}$ shows that it differs from 
$\br{\diag{3mm}{2}{2}{
  \piccirclearc{1 -1}{1.41}{45 135}
  \piccirclearc{1 3}{1.41}{-135 -45}
}}$ just by multiplication and addition of a unit, hence the
determinants of $L_1=\diag{3mm}{2}{2}{
  \piccirclearc{-1 1}{1.41}{-45 45}
  \piccirclearc{3 1}{1.41}{135 -135}
}$ and $L_2=\diag{3mm}{2}{2}{
  \piccirclearc{1 -1}{1.41}{45 135}
  \piccirclearc{1 3}{1.41}{-135 -45}
}$ differ just by $\pm 1$. Then
these determinants are $n$ and $n+1$ respectively (the even one
corresponding to the 2 component link and the odd one to the knot).
To see this, either evaluate \eqref{tlV} and use the skein
relation \eqref{1} for $V$ for the unknotting crossing change
at $t=-1$, or alternatively use the relation \eqref{Qrel} for $Q$
evaluated at $z=2$ and the property of \cite{BLM} that $\dt(K)=\sqrt
{Q(2)}$. Then Consider the $Q$ relation at $z=\frac{\sqrt{5}-1}{2}$ at
the unknotting crossing, and let $k$ and $l$ be the two numbers
$\dim H_1(D_{L_{1,2}},\bZ_5)$. The claim follows using the result of
\cite{Jones2} and Rong \cite{Rong} on the form of
$Q\bigl(\frac{\pm\sqrt{5}-1}{2}\bigr)$ by
rewriting the resulting condition as a congruence modulo the minimal
polynomial of $\frac{\pm\sqrt{5}-1}{2}$. \qed

The criterion for $Q$ can be generalized to higher unknotting numbers.

\begin{theorem}\label{T3.1}
If $Q_K\bigl(\frac{\sqrt{5}-1}{2}\bigr)=-\bigl(-\sqrt{5}\bigr)^k$,
then $u(K)>k$.
\end{theorem}

\proof The proof goes through by slight modification of the arguments
in the proof of Proposition \reference{th3}, and via induction on the
unknotting number, as in \cite{Traczyk}. One needs to observe,
using the relation of $Q$, that if $K$ and $K'$ differ by a crossing
change, then the value $Q_K\bigl(\frac{\sqrt{5}-1}{2}\bigr)/
Q_{K'}\bigl(\frac{\sqrt{5}-1}{2}\bigr)\in\{\pm 1,-5^{\pm 1/2}\}$, but
it cannot be $+5^{\pm 1/2}$. 
\qed


\section{The Goeritz matrix and linking form, and a
criterion for unknotting number 3}

Let us turn to a formerly known
topological approach using the linking form $\lm$ on $H_1(D_K)$. It is
a consequence of a result of Montesinos \cite{Montesinos}
and was first stated by Lickorish \cite{Lickorish}. 
(Note that $D_K$ inherits an orientation from $S^3$, so that
the linking form is given by $+U^{-1}$, where $U$ is a Goeritz
\cite{GorLit} matrix for $K$.)

\begin{prop} (Montesinos-Lickorish)\label{TML}
If $u(K)=1$, and $\dt(K)=D$ then ($H_1(D_K)$ is cyclic
and) there is a generator $g$ of $H_1(D_K)$, with $\lm(g,g)=\pm 2/D
\in\bQ/\bZ$.
\end{prop}

Here $\lm$ denotes the linking form on $H_1(D_K)$ and $u_{\pm}$ denotes
the \em{signed} unknotting number (see \cite{Traczyk,CL}), that is,
$u_+(K)=1$ iff $K$ unknots by switching a positive crossing to the
negative. This is the property, to which we found obstructions in
\S\reference{JP}, and we will use these obstructions shortly.

There is an extension of the linking form criterion which may
be useful when considering signed unknotting number. It appears
in a similar form to the one given here in \cite[proposition 2.1]
{CL}. (Compare also the proof in Kawauchi's book 
\cite[theorem 11.2.3, pp.\ 147-8]{Kawauchi}.)

\begin{theorem}\label{rML}
Let $K$ be an unknotting number one knot which can be unknotted by
switching a positive crossing to the negative, and $\dt(K)=D$.
Then there is a generator $g$ of
$H_1(D_K)$ with $\lm(g,g)=+2/D\in\bQ/\bZ$ if $\sg(K)=0$ and
$\lm(g,g)=-2/D$ if $\sg(K)=2$.
\end{theorem}

\proof $K$ and the unknot possess Seifert matrices $S$ and $S'$ with
$S_{i,j}=S'_{i,j}$ for $(i,j)\ne (1,1)$ and $S_{1,1}=S'_{1,1}+1$.
Now $S+S^T$ is a representation matrix for $H_1(D_K)$, and differs
just in the upper left entry from $S'+S'^T$. The determinant $\dt(
S+S^T)$ is up to sign equal to $\dt(K)=D$, and $\dt(S'+S'^T)=\pm 1$.
Thus the minor of $S+S^T$ obtained by deleting the first row and column
from $S+S^T$ has determinant $(\pm D\pm 1)/2$. The homology
element $g'$ represented by the first row and column then has
linking form $\lm(g',g')=(\pm D\pm 1)/2D$. Set $g=2g'$. Then
$\lm(g,g)=\pm 2/D$, and hence $g$ is a generator. The important point
to notice now is that the sign of $\lm(g,g)$ depends on whether 
$\dt(S+S^T)$ and $\dt(S'+S'^T)$ have the same sign or not, and that
this is equivalent to the signature condition. \qed

\begin{remark}\label{rK'}
More generally the proof shows that if $K$ can be turned into some
knot $K'$ by switching a positive crossing to negative and
$\dt(K')=D'$, then there is a (not necessarily generating)
element $g\in H_1(D_K)$ with $\lm(g,g)=+2D'/D$ if $\sg(K)=
\sg(K')$ and $\lm(g,g)=-2D'/D$ if $\sg(K)=2+\sg(K')$. (This remark
will be used later in the proof of Theorem \reference{tH2}.)
\end{remark}

\begin{rem}\label{XX}
In \cite{CL} an example is given, the $(-9,5,-9)$-pretzel knot,
on which Proposition \reference{TML} does not apply, but theorem
\reference{rML} excludes one of the possible signs of an
unknotting crossing.
\end{rem}

%
%
%


The (conjectured) coincidence of the criteria with $Q$ and $\lm$ for
unknotting number $1$ (for determinant divisible by $5$) suggests
that further relations between both may exist also for higher
unknotting numbers. And indeed, Theorem \reference{rML}, although 
apparently in practice not much more effective than its unsigned
version, can in theory be pushed further, at least in a special case,
to prove a new criterion, this time for unknotting number 3.

\begin{theorem}\label{tH2}
If $d=\dt(K)$ is a square and has no divisors of the
form $4k+3$, and $\sg(K)=4$, then $u(K)>2$.
\end{theorem}

%

\proof
If $u(K)=2$, then there exists an
unknotting sequence $K\to K'\to\bigcirc$ (latter denoting the unknot).
It is known, originally from \cite{Murasugi}, that for any knot $K$
\[
\dt(K)\,
\left\{\begin{array}{c}
\equiv 1\,(4) \\
\equiv 3\,(4) \end{array}\right.\mbox{\quad iff\quad}
\sg(K)\,
\left\{\begin{array}{c}
\equiv 0\,(4) \\
\equiv 2\,(4) \end{array}\right.\,.
\]
Clearly $\sg(K')=2$, and hence $d'=\dt(K')\equiv 3\bmod 4$.
But then $K$ and $\bigcirc$ are obtained from $K'$ by
switching crossings of opposite sign (under both of which the
signature changes). Now, we apply the argument proving
Theorem \reference{rML} for the crossing switch $K'\to\bigcirc$,
and the modified version of this argument given in Remark
\reference{rK'} to the crossing switch $K'\to K$ (where in fact we
did not involve anything more than the determinant of the Seifert
matrices). By combining both, we find some $c\in\bZ_{d'}$ with $c^2
\equiv -d\bmod d'$. This is impossible, however, whatever $d'$ may
be, under the condition on $d$ we assumed. (Consider the
congruence just modulo a prime $4k+3$ dividing $d'$ and use the
First supplementary law of number theory for Legendre's symbol.) \qed

\begin{remark}\label{RK'}
As in remark \reference{rK'}, the same argument shows that
more generally, if $\sg(K)-\sg(K')=4$ and the product $\dt(K)\cdot
\dt(K')$ is a square and has no divisors of the form $4k+3$
(which in particular means that $4\mid \sg(K),\,\sg(K')$),
then $K$ and $K'$ cannot be interconverted by 2 crossing changes
(that is, have distance at least 3 in the sense of Murakami
\cite{Murakami}). For example, $5_1$ cannot be made into
$4_1$ by two crossing changes.
\end{remark}

\section{Applications}
 
If the criterion in Proposition \reference{thV} applies to both
mirror images, or considering one of them is irrelevant (e.g.
because the Jones polynomial is self-conjugate or because of the
signature), it successfully, and easily, excludes unknotting number
one, as for the following knots (most of them due to Nakanishi and
Kanenobu--Murakami), thus shown to have unknotting number two: $7_4$,
$8_{18}$, $9_{15}$, $9_{17}$, $9_{37}$, $9_{40}$, $9_{46}$, $9_{47}$
and $9_{48}$. Even for unknotting number one knots, as $6_1$ and $7_7$,
the criterion can give non-trivial information, namely, that such knots
are unknottable by switching crossings of only one sign
(dependingly on how they are obversed), despite having $\sg=0$. 
Finally, the criterion in Proposition \reference{thV} recovers
after \cite{KanMur} Lickorish's result \cite{Lickorish} for $7_4$, and
also Traczyk's recent example $10_{67}$ \cite{Traczyk}.

\begin{rem}
K.\ Murasugi informed me of a paper of Miyazawa \cite{Miyazawa},
where he obtained the unknotting numbers of some 10 crossing knots
using a formula relating the value $V'(-1)$, determinant and the Conway
degree-4-Vassiliev invariant (Theorem 8). (Their entries were completed
in Kawauchi's tables.) For them all other criteria~-- Lickorish's,
Traczyk's and ours, also apply. Thus it appears that Traczyk's
criterion for unknotting number one (and so basically ours, too) is
equivalent to Miyazawa's (Miyazawa also obtains a slightly weaker
version of corollary \reference{cMi} in his corollary 7), and they
are implied by the linking form condition of Lickorish (see conjectures
in \S\reference{Scj}).
\end{rem}

The results 
we proved in \S\reference{SQ} allow to settle the unknotting numbers
of the knots in \eqref{teq} as follows: $9_{49}$ and $10_{103}$,
and also $5_1\#4_1$, have $u=3$ by Theorem \reference{T3.1}, and
$8_{16}$, $10_{86}$, $10_{106}$, $10_{109}$, $10_{116}$,
and $10_{121}$ have $u=2$ both by Theorem \reference{T3.1}
and Proposition \reference{th3}.

As the simplest special case of this proposition, any knot $K$ with 
$Q_K\bigl(\frac{\sqrt{5}-1}{2}\bigr)=+\sqrt{5}$ cannot have unknotting
number one. Another possible special case would be a knot $K$ with $\det
(K)\equiv \pm 2 (5)$ and $Q_K\bigl(\frac{\sqrt{5}-1}{2}\bigr)=+1$
respectively 
$\dt(K)\equiv \pm 1 (5)$ and $Q_K\bigl(\frac{\sqrt{5}-1}{2}\bigr)=-1$.
However, it was shown by Rong \cite{Rong} that such knots do not exist,
that is, the sign of $Q_K\bigl(\frac{\sqrt{5}-1}{2}\bigr)$ for $5\nmid
\dt(K)$ is exactly the same as $Q_K(2) \bmod 5$. (Rong gives an
interpretation of the sign in terms of the Goeritz form.) Therefore,
Proposition \reference{th3} is basically equivalent
to Theorem \reference{T3.1} for $k=1$.

The unknotting numbers (equal to 2) of 
$10_{86}$, $10_{105}$, $10_{106}$, $10_{109}$,
$10_{116}$ and $10_{121}$ are applications of
proposition \reference{TML}. Lickorish's proof of $u(7_4)=2$ also
consists in application of this condition, as well as this for
$8_{16}$, a result of J.\ R.\ Rickard whose proof was never published.
The disadvantage of this method is that the calculation by hand of $\lm$ is
in general not pleasant. Nevertheless, a nice approach for calculating
$\lm$ was developed by Gordon and Litherland \cite{GorLit} via the
Goeritz matrix $U$. It was carried out on the example $8_{16}$ by
Murakami and Yasuhara in \cite{MurYas} (in slightly generalized form)
to give $\lm(g',g')=\pm 11/35$ for some generator $g'$, recovering
Rickard's result (as $\pm 11$ is not twice a square modulo $35$). From
the Goeritz matrix, $\lm$ is given by the matrix $U^{-1}$ on the
generating set of $H_1(D_K)$, on which $U$ gives the relation matrix.

\begin{exam}
To make the calculations verifiable, we give as example the Goeritz
matrix, the Dowker--Thistle\-thwaite notation, signature, and a
generator of $H_1(D_K)$ (in the basis determined by the Goeritz matrix)
for the knot $K=10_{105}$ obtained from its (unique) 10 crossing
diagram.

{
\begin{verbatim}
  5  0 -2 -2
  0  4 -2 -1
 -2 -2  5  0
 -2 -1  0  3
 10      105 :   91 (0,0,0,1)
 10      105      4  12  16  20  18   2   8   6  10  14        91    2
\end{verbatim}
}
\end{exam}

\begin{rem}
In Rolfsen's tables $10_{83}$ and
$10_{86}$ are swapped: the Conway notation and Alexander polynomial
for each one refers to the diagram of the other. The convention for
$10_{86}$ here is that the Conway notations and Alexander polynomials
are swapped to fix the discrepancy, and \em{not} the diagrams, as in
\cite{Kawauchi} (so our $10_{86}$ is Kawauchi's $10_{83}$).
\end{rem}

As mentioned, theorem \reference{rML} is in general not
much more powerful than proposition \reference{TML}.
However, its consequence, Theorem \reference{tH2} again applies
to $9_{49}$ and $5_1\#4_1$ (but not to $10_{103}$) to show $u=3$.
Strangely, for all prime knots of $\le
16$ crossings with $\sg=4$, $d=25$ or $d=625$ and 2 torsion numbers
in $H_1(D_K)$, we had $Q\bigl(\frac{\sqrt{5}-1}{2}\bigr)=-5$, so that
the $Q$ method excludes $u=2$ as well. This gives
further hints to an intimate relationship between the sign of the
Jones-Rong value of $Q$ and the topology of the double cover
of the knot complement. 

Nevertheless, the criterion of Theorem \reference{tH2}
will clearly give new information for other primes $4k+1$, which
cannot be captured by $Q$. For example, the knot $12_{664}$ of the
tables available in KnotScape \cite{%
KnotScape} has $H_1=\bZ_{13}\oplus \bZ_{13}$ and $\sg=4$, and hence
$u=3$, which cannot be shown by any other method I know of. Another
such, possibly more common, example is $6_3\# 7_3$.


As a summary, beside $8_{16}$ and $9_{49}$, there are several
undecided 10 crossing knots in Kawauchi's tables \cite{Kawauchi},
whose unknotting number we found. The (most likely correct) list
of prime knots of at most 9 crossings with unknown unknotting number
now is: $8_{10}$, $9_{29}$ and $9_{32}$ with $u\in\{1,2\}$ and
$9_{10}$, $9_{13}$, $9_{35}$, $9_{38}$ with $u\in\{2,3\}$
(see \cite[table, p.~49]{Kirby}).
Beside the knots so far mentioned, we can complete the entry for
$10_{131}$, whose unknotting number is 1, as show both its diagrams
in KnotScape and in Rolfsen's book.
%
The unknotting numbers (so far known to me) of the
simplest composite knots are compiled in a table in the appendix.
 
\section{Conjectures and problems\label{Scj}}



In this final section we summarize the unexplained phenomena that
came up in our attempts to find unknotting number criteria studying
the linking form and the polynomial values, and most of which
have been implicitly suggested in our previous discussion.
(However, there are also some new ones.)
All they are supported by strong empirical evidence.
Some of these statements indicate, that the conditions for
the unknotting number given by Lickorish \cite{Lickorish}, Traczyk
\cite{Traczyk}, and above in this paper seem to imply others,
at least in special cases. 

The first problem came up in the study of the structure of $H_1$
of knots to which Theorem \reference{tH2} is applicable.

\begin{conjecture}
There are no knots with $\sg=4$ and cyclic $H_1$ of order a
prime square (i.e., $H_1$ is always a double in all such cases).
If $\sg(K)=4$ and $H_1=\bZ_5\oplus\bZ_5$, then
$Q_K\bigl(\frac{\sqrt{5}-1}{2}\bigr)=-5$ (rather than $+5$).
\end{conjecture}

At least the first part of this conjecture
is not true for non-prime squares, even for higher
(even) prime powers~-- there are  for example two 16 crossing knots
with $\sg=4$ and $H_1=\bZ_{25}\oplus\bZ_5\oplus\bZ_5$. By
\cite[theorem 3.10]{HNK} there is no $\sg=4$ knot with determinant 1.
As we saw, the simplest examples supporting the second
part of this conjecture are $9_{49}$ and $5_1\#4_1$.

\begin{conjecture} 
If $K$ has prime determinant $\dt(K)=D$ and signature 0,
then there is an element $g\in H_1$ with $\lm(g,g)=\pm 2/D$.
\end{conjecture}

For $\sg=2$ this is true, as then $D=4k+3$. If such $D$
is prime, the multiplicative
group $\bZ_D^*$ is cyclic and has no square roots of $-1$, so that
each residue class $\bmod D$ is either of the form $2a^2$ or $-2a^2$.
For $D=4k+1$ one half of the residue classes are of both forms, and
the other half of none of them. Latter seem never to occur for $\lm$
if $\sg=0$.


\begin{conjecture} 
If $K$ has cyclic $H_1$ of order $\dt(K)=D$ divisible by $5$, then
\[
Q_K\Bigl(\frac{\sqrt{5}-1}{2}\Bigr)\,=\,
\left\{\begin{array}{cc}
-\sqrt{5} & \mbox{if\ }\mathrel{\,\exists}\ g\in H_1\,:\,\lm(g,g)=\pm 2/D \\
+\sqrt{5} & \mbox{if\ }\not\mathrel{\,\exists}\ g\in H_1\,:\,\lm(g,g)=\pm 2/D 
\end{array}\right.\,.
\]
\end{conjecture}

\begin{conjecture}\label{cLT} 
Assume $\sg(K)=+2$, and $H_1$ is cyclic. (a) If $V_K\left(e^{\pi i/3}
\right)=-i\sqrt{3}$, then there is no $g\in H_1$ with $\lm(g,g)=\pm
2/D$. (b) If there is no $g\in H_1$ with $\lm(g,g)=-2/D$, then
there is neither a $g\in H_1$ with $\lm(g,g)=+2/D$.
\end{conjecture}

It seems possible that if $\sg(K)=+2$ and $H_1$ is cyclic
of order divisible by $3$, then the conditions
$V_K\left(e^{\pi i/3}\right)=-i\sqrt{3}$ and that
there is no $g\in H_1$ with $\lm(g,g)=\pm 2/D$ are equivalent.
This motivates the first part of the conjecture.
In case $\sg=2$, I know of no example for which the
signed unknotting number information of Theorem \reference{rML}
contradicts that of the signature, but Proposition \reference{TML}
does not apply. This is the origin of the second part
(for the case $\sg=0$, see remark \reference{XX}).

It seems difficult to generalize Theorem \reference{tH2} in some way.
The condition $\sg=4$ (rather than $\sg=0$) and the lack of divisors
$4k+3$ of the determinant are both necessary (otherwise consider the
connected sums of twist knots). How to weaken the squareness condition
on the determinant is not clear either, since a computer experiment
revealed that all non-square numbers $4k+1<400$ are realized as
determinants of knots with $u=2$ and $\sg=4$.

The following is a problem on a class of knots where the
unknotting number one condition has been strengthened.

\begin{question}
The (reduced alternating)
trefoil and figure eight knot diagrams have the property that they
unknot by switching \em{any arbitrary} crossing. Are they the
only (non-trivial\footnote{Diagrams of the unknot with all crossings
reducible trivially have this property, too.}) knot diagrams with
this property? More generally, are for $k>1$ the $(2,2k\pm 1)$-torus
knot diagrams (with one possible kink in the `$-$' case) the only
diagrams which unknot by switching any arbitrary collection of $k$
crossings?
\end{question}

So far the only observation towards the case $k=1$ is that except for
the trefoil diagram any other such diagram must have $\sg=0$ and
$|V(-1)|\equiv 1$ or $5\bmod 12$.

As a final remark on the $V$ criterion, we mention that there are still
some possibilities left open. The most promising way appears to be to
consider the Homfly polynomial of the 2-component link arising by
smoothing the unknotting
crossing and to attack the existence of such a polynomial by the
various Vassiliev invariant identities worked out by Kanenobu.
We may record an interesting outcome of this (or some
similar) idea at a later stage.

Also, one may try to find more special evaluations of the polynomials
related to branched cover homology, but the results of \cite{Rong} and
\cite{qeval} for the $Q$ polynomial (which hold in an analogous
form also for $V$, see \cite{JVW}) suggest, that
such evaluations, beside the known ones, are very unlikely to exist.

\noindent{\bf Acknowledgement.} The investigation of this
paper was inspired by Traczyk's paper \cite{Traczyk} and
his talk on the Conference on  Knot Theory ``Knots in Hellas,
98''. I would also wish to thank to T.~Cochran, W.~B.~R.~Lickorish and
H.~Murakami for helpful discussions on the linking form and to
K.\ Murasugi for telling me of the paper \cite{Miyazawa}.

\begin{appendix}

\renewcommand{\thesection}{\Roman{section}}
\def\@seccntformat#1{Appendix \csname the#1\endcsname . \quad}

\section{Tables}

The unknotting number is not \em{a priori} additive
under connected sum. (That this is true is a long--standing conjecture
proved only in the first non-trivial case by Scharlemann \cite{Scharl}
and possibly even false in general.) Therefore, there seems no reason to favorize
prime knots in the tabulation of unknotting numbers. For this reason I include a
table of the unknotting numbers of composite $\le 10$ crossing knots
as far as feasible to me.

I adopt the convention of Traczyk \cite{Traczyk} that the number
$x$ followed by $y$ copies of `?' means that the unknotting number is
at most $x$, and very likely (or in the case $y=0$ known to be) $x$,
but that the values $x-y,\dots, x-1$ have not yet been excluded.
The knots are recorded up to mirroring, taking Thistlethwaite's
obversion convention for the factors (i.e., if one of $K$ and $!K$
is positive, then always $K$ is taken to be such rather than $!K$).

To compile these values, it basically suffices to apply the standard
results of \cite{Murasugi,Scharl,Wendt}. These methods work except for
three knots, where other methods are needed: $7_4\#3_1$ and
$3_1\#3_1\#4_1$ (Traczyk) and $4_1\# 5_1$ (see above).

\[
\begin{mytab}{|l|l||l|l||l|l||l|l|}{ & & & & & & & }
\hline[1mm]%
$K$ & \multicolumn{1}{|c||}{$u(K)$} & $K$ & \multicolumn{1}{|c||}{$u(K)$}
& $K$ & \multicolumn{1}{|c||}{$u(K)$} & $K$ &
\multicolumn{1}{|c|}{$u(K)$} \\[2mm]%
\hline
\hline[2mm]%
$3_1\#3_1$ & 2	& $3_1\#!6_1$ & 2  & $3_1\#!7_2$ & 2 & $3_1\#!7_7$ & 2 \\
$3_1\#!3_1$ & 2 & $3_1\#6_2$ & 2   & $3_1\#7_3$ & 3 & $4_1\#6_1$ & 2 \\
$3_1\#4_1$ & 2	& $3_1\#!6_2$ & 2  & $3_1\#!7_3$ & 3? & $4_1\#6_2$ & 2 \\
$3_1\#5_1$ & 3	& $3_1\#6_3$ & 2   & $3_1\#7_4$ & 3 & $4_1\#6_3$ & 2 \\
$3_1\#!5_1$ & 3?& $4_1\#5_1$ & 3   & $3_1\#!7_4$ & 3? & $5_1\#5_1$ & 4 \\
$3_1\#5_2$ & 2  & $4_1\#5_2$ & 2   & $3_1\#7_5$ & 3 & $5_1\#!5_1$ & 4?? \\
$3_1\#!5_2$ & 2 & $3_1\#3_1\#4_1$ & 3 & $3_1\#!7_5$ & 3? & $5_1\#5_2$ & 3 \\
$4_1\#4_1$ & 2	& $3_1\#!3_1\#4_1$ & 3? & $3_1\#7_6$ & 2 & $5_1\#!5_2$ & 3? \\
$3_1\#3_1\#3_1$ & 3 & $3_1\#7_1$ & 4 & $3_1\#!7_6$ & 2 & $5_2\#5_2$ & 2 \\
$3_1\#3_1\#!3_1$ & 3 & $3_1\#!7_1$ & 4?? & $3_1\#7_7$ & 2 & $5_2\#!5_2$ & 2 \\
$3_1\#6_1$ & 2   & $3_1\#7_2$ & 2 &              &   &              & \\
\hline
\end{mytab}
\]

\end{appendix}

\end{document}